\documentstyle[amssymb,amsfonts,12pt]{amsart}

\newtheorem{theorem}{Theorem}[section]
\newtheorem{lemma}[theorem]{Lemma}
\newtheorem{corollary}[theorem]{Corollary}
\newtheorem{definition}[theorem]{Definition}

\def\C{{\mbox{\rm\kern.24em
\vrule width.03em height1.43ex depth-.052ex \kern-.26em C}}}
\def\QSet{\mbox{\rm\kern.24em
\vrule width.03em height1.48ex depth-.051ex \kern-.26em Q}}
\def\Z{{\bf Z}}
\def\R{{\mbox{\rm I\kern-.22em R}}}
\def\P{{\bf P}}
\def\Q{{\bf Q}}
\def\T{{\bf T}}

\def\size{{\rm size}}
\def\diam{{\rm diam}}
\def\sgn{{\rm sgn}}
\def\M{{\rm M}}

\def\eps{\varepsilon}
\def\pv{{\vec P}}
\def\Pv{{\vec{\bf P}}}
\def\dist{{\rm dist}}

\def\111{\gamma}
\def\hyperplane{\Gamma}

\def\subspace{{\Gamma'}}
\def\orthogonal{{\Gamma''}}

\def\be#1{\begin{equation}\label{#1}}
\def\bas{\begin{align*}}
\def\eas{\end{align*}}
\def\bi{\begin{itemize}}
\def\ei{\end{itemize}}
\newenvironment{proof}{\noindent {\bf Proof} }{\endprf\par}
\def \endprf{\hfill  {\vrule height6pt width6pt depth0pt}\medskip}
\def\emph#1{{\it #1}}

\title{Multi-linear operators given by singular multipliers}

\author{Camil Muscalu}
\address{Department of Mathematics, Brown University, Providence RI 02912}
\email{camil@@math.brown.edu}

\author{Terence Tao}
\address{School of Mathematics, UNSW, Sydney NSW 2052 AUSTRALIA}
\email{tao@@math.ucla.edu}

\author{Christoph Thiele}
\address{Department of Mathematics, UCLA, Los Angeles CA 90095-1555}
\email{thiele@@math.ucla.edu}

\begin{document}

\begin{abstract}  We prove $L^p$ estimates for a large class of multi-linear operators, which includes the multi-linear paraproducts studied by Coifman and Meyer \cite{coifmanm6}, as well as the bilinear Hilbert transform.
\end{abstract}

\maketitle

\section{Introduction}

Let $n > 1$ be an integer, and let $m(\xi_1, \ldots, \xi_n)$ be a function on the $n-1$-dimensional vector space
$$ \hyperplane = \{ \xi \in \R^n: \xi_1 + \ldots + \xi_n = 0 \}.$$
For any $m$, we associate the multi-linear operator $T = T_m$ on $n-1$ functions
on $\R$ by
\be{tm-def}
 T_m(f_1, \ldots, f_{n-1}\hat{)}(-\xi_n) = \int \delta(\xi_1 + \ldots + \xi_n) m(\xi) \hat f_1(\xi_1) \ldots \hat f_{n-1}(\xi_{n-1})\ d\xi_1 \ldots d\xi_{n-1},
\end{equation}
where $\xi = (\xi_1, \ldots, \xi_n)$.  We may write this operator more symmetrically as an $n$-linear form $\Lambda = \Lambda_m$ given by
$$ \Lambda_m(f_1, \ldots, f_n) = \int \delta(\xi_1 + \ldots + \xi_n)
m(\xi) \hat f_1(\xi_1) \ldots \hat f_n(\xi_n)\ d\xi;$$
the relationship between $T$ and $\Lambda$ is given by
\be{lambda-t} \Lambda(f_1, \ldots, f_n) = \int T(f_1, \ldots, f_{n-1})(x) f_n(x)\ dx.
\end{equation}

When $n=2$, $T$ is a Fourier multiplier, and it is well known that such operators are bounded on $L^p$, $1 < p < \infty$, if $m$ is a symbol of order 0.
Coifman and Meyer \cite{coifmanm1}-\cite{coifmanm6},
Kenig and Stein \cite{kenigs}, and Grafakos and Torres \cite{grafakost} 
extended this result to the 
$n > 2$ case, showing that
one had the mapping properties
\be{holder} T: L^{p_1} \times \ldots \times L^{p_{n-1}} \to L^{p'_n}
\end{equation}
whenever 
\begin{equation}\label{banach}
1 < p_i \le \infty
\end{equation}
 for $i=1, \ldots, n-1$, 
\begin{equation}\label{output}
1/(n-1)<p_n'<\infty\ \ \ ,
\end{equation}
and
\be{scaling} \frac{1}{p_1} + \ldots + \frac{1}{p_n} = 1,
\end{equation}
and $m$ satisfies the symbol estimates
\be{symb}
|\partial_\xi^\alpha m(\xi)| \lesssim |\xi|^{-|\alpha|}
\end{equation}
for all partial derivatives $\partial_\xi^\alpha$ on $\Gamma$ up to some 
finite order.  
We interpret estimate (\ref{holder}) in the way that $T$ is originally
defined on the product of suitable subspaces of the $L^{p_i}$
and then extends to the product of the closures of these subspaces. 
In case ${p_i}\neq \infty$ the subspace is simply the test function 
space which is dense in $L^{p_i}$.  If any estimate of the type (\ref{holder})
holds with $p_i\neq \infty$ for all $i$, then we can use this to 
unambiguously define $T$ on the product of $n-1$ copies
of $L^1\cap L^\infty$. Once this is done, we can choose
$\overline{L^1\cap L^\infty}$ as subspace of $L^\infty$ whenever
$p_i=\infty$ in some other eponent tuple.
If ${p_n}'\ge 1$, we can use a duality argument to extend the operator from 
$\overline{L^1\cap L^\infty}$ to $L^\infty$.
If ${p_n}'<1$ we shall be satisfied with replacing
$L^\infty$ by $\overline{L^1\cap L^\infty}$ in \eqref{holder}
where applicable.

The interesting observation that $p_n'$ can 
be smaller or equal $1$ traces back (at least) to papers by C. Calderon 
\cite{calderon} and Coifman and Meyer \cite{coifmanm1}, where special
multilinear operators are discussed.

When $m$ is identically one then $T_m$ is the pointwise product operator
$$ T(f_1, \ldots, f_{n-1}) = f_1 \ldots f_{n-1}\ \ \ ,$$
so estimate \eqref{holder} may be viewed as a 
generalization of H\"older's inequality, 
where products are replaced by paraproducts.

The bilinear Hilbert transform
$$ T(f_1,f_2) = \int f_1(x-t) f_2(x+t) \frac{dt}{t}$$
can also be viewed as an operator of the form \eqref{tm-def}, with symbol
$$ m(\xi_1,\xi_2,\xi_3) = \pi i \sgn(\xi_2 - \xi_1)\ \ \ .$$
This multiplier does not satisfy the estimates \eqref{symb}.  Nevertheless, Lacey and Thiele \cite{laceyt1}, \cite{laceyt2} showed that \eqref{holder} continues to hold, provided
that one makes the additional assumption $p'_3 > 2/3$.

The purpose of this paper is to unify these results, allowing us to prove \eqref{holder} for a class of multipliers which are singular on a subspace of $\hyperplane$.  More precisely, we have

\begin{theorem}\label{main}  Let $\subspace$ be a subspace of $\hyperplane$ of dimension $k$ where
\be{k-constraint}
0 \leq k < n/2. 
\end{equation}
Assume that $\subspace$ is non-degenerate in the sense that for every
$1 \leq i_1 < i_2 < \ldots i_k \leq n$, the space $\subspace$ is a graph
over the variables $\xi_{i_1}, \ldots, \xi_{i_k}$.
Suppose that $m$ satisfies the estimates
\be{symb-gamma}
|\partial_\xi^\alpha m(\xi)| \lesssim \dist(\xi,\subspace)^{-|\alpha|}
\end{equation}
for all partial derivatives $\partial_\xi^\alpha$ on $\Gamma$ up to some finite order.  Then \eqref{holder} holds whenever
\eqref{banach}, \eqref{output}, \eqref{scaling} hold and
\be{odd}
\frac{1}{p_{i_1}} + \ldots + \frac{1}{p_{i_r}} < \frac{n-2k+r}{2}
\end{equation}
for all $1 \leq i_1 < \ldots < i_r \leq n$ and $1 \leq r \leq n$.
\end{theorem}

In particular, \eqref{holder} holds whenever 
$1 < p_i \le \infty$ for $i = 1, \ldots, n$ and \eqref{scaling} holds.

As discussed above the case $k=0$ is well known.
The Lacey-Thiele theorem is covered by the case $n=3$, $k=1$.  Unfortunately this theorem does not quite cover the trilinear Hilbert transform
$$ T(f_1,f_2,f_3) = \int f_1(x-t) f_2(x+t) f_3(x+2t) \frac{dt}{t}$$
since one has $n=4$, $k=2$ in this case, which does not satisfy \eqref{k-constraint}.  To obtain the analogue of this theorem when \eqref{k-constraint} fails would probably require radically different techniques than the ones developed to date.  However, an elementary argument can be used to handle this case if enough functions are in the Wiener algebra $A = {\cal F}^{-1} L^1$; see 
Section \ref{remarks}.

In the $k=0$ case the origin $\xi = 0$ has special significance, and this is reflected in the tools used to handle this case, namely Littlewood-Paley theory and/or wavelets.  However, when $k > 0$ there is no preferred frequency origin,
and the tools used should be invariant under frequency translations along $\subspace$.  This necessitates the employment of ``tiles'' in the time-frequency plane which have arbitrary frequency location, spatial location,
and scale.

If the multiplier $m$ of Theorem \ref{main}
is invariant under translations in direction of $\subspace$, then we 
can write the $n-1$-linear operator $T$ as
$$
T(f_1,\dots ,f_{n-1})(x)=\int_{\orthogonal \cap \hyperplane}
f_1(x+\gamma_1)\dots f_{n-1}(x+\gamma_{n-1})K(\gamma)
\, d\gamma\ \ \ ,
$$
where $\orthogonal$ is the orthogonal complement of $\subspace$, 
$d\gamma$ is Lebesgue measure
on $\orthogonal\cap \hyperplane$,  $\gamma_i$ is the $i$-th coordinate
of $\gamma$ as element of $\R^n$, and $K$ is a Calderon-Zygmund 
kernel on the space $\orthogonal\cap \hyperplane$.
Thus we obtain $L^p$ bounds for such operators provided 
$n-1\le 2d$, where $d$ is the dimension of $\orthogonal\cap \hyperplane$.
This gives a partial answer
to question (2) in \cite{kenigs} raised by Kenig and Stein.

It would be interesting to study the behaviour of the bounds
in \eqref{holder} as the space $\subspace$ degenerates in the sense
of Theorem \ref{main}, see \cite{thiele} for some results in this direction
in the special case of the bilinear Hilbert transform. We do not discuss 
this issue here.

This paper is organized as follows.  In Section \ref{interp-sec} we introduce some multi-linear interpolation theory, which allows us to reduce \eqref{holder} to a ``restricted type'' estimate on the $n$-form $\Lambda$.  In Section \ref{except} we remove an exceptional set, and reduce matters to estimating $\Lambda$ on functions whose Hardy-Littlewood maximal function is under control.  In Section \ref{discrete} we then decompose the multiplier $m$ using a Whitney decomposition, which allows us to replace $\Lambda$ by a discretized analogue which involves the size of the $f_i$ on various tiles in the time-frequency plane; roughly speaking, we only need consider those tiles that lie outside the exceptional set.  
To handle these tiles we first consider the case $k=1$.
This is done by subdividing the tiles into essentially disjoint
trees, using Littlewood-Paley theory to estimate the contribution of each
tree individually, and then using orthogonality arguments to control the
total number of trees.  Finally, in Section \ref{induct-sec}, we induct on $k$
to obtain the general case.

The first author wishes to express his gratitude to the UCLA
Department for its hospitality during his visit to Los Angeles
and to Jill Pipher for her financial and moral 
support.  The second author thanks Michael Lacey and Jim Wright for many helpful discussions during a delightful semester at the University of 
New South Wales.
The second and third authors are supported by NSF grants 
DMS 9706764 and DMS 9970469 respectively. The third author
acknowledges an enjoyable stay at the Erwin Schr\"odinger Institute
in Vienna, during which part of this work was done.

\section{Preliminaries}

We use $A \lesssim B$ to denote the statement that $A \leq CB$ for some
large constant $C$, and $A \ll B$ to denote the statement that $A \leq C^{-1} B$
for some large constant $C$.

Throughout the paper we shall assume $k>0$, the case $k=0$ being well known.

We make the a priori assumption that the symbol $m$ is smooth and compactly
supported; this makes $\Lambda$ bounded on every product of $n$
Lebesgue spaces.  Our final estimates will not depend on the smoothness or support bounds that $m$ satisfies, and the general case can be handled by
the usual limiting argument.

If $I$ is an interval, then $CI$ denotes the interval with the same center but $C$ times the length.  Let $\chi_I$ denote the characteristic function of $I$.
We define the approximate cutoff function $\tilde \chi_I$ as
$$ \tilde \chi_I(x) = (1 + \frac{\dist(x,I)}{|I|})^{-1}.$$

We use $\M f$ to denote the Hardy-Littlewood maximal function.

\section{Interpolation}\label{interp-sec}

In this section we develop some multi-linear interpolation theory which allows us to reduce \eqref{holder} to a certain ``restricted type'' estimate 
on $\Lambda$.

In this section we find it convenient to work with the 
quantity $\alpha_i=1/p_i$ when $p_i$ is the exponent of $L^{p_i}$.
We fix $n$ throughout this section.

\begin{definition}  A tuple $(\alpha_1,\dots ,\alpha_n)$ is called
admissible, if 
\begin{equation}\label{admissible}
-\infty<\alpha_i<1
\end{equation}
for all $1\le i\le n$,
\begin{equation}\label{admissum}
\sum_i\alpha_i=1\ \ \ ,
\end{equation}
and there is at most one index $j$ such that $\alpha_j<0$. 
We call an index $i$ {\it good} if $\alpha_i\ge 0$, and we call it
{\it bad} if $\alpha_i<0$. A {\it good} tuple is an admissible 
tuple without bad index, a {\it bad} tuple is an admissible tuple with 
a bad index.
\end{definition}

\begin{definition}  Let $E$, $E'$ be sets of finite measure.
We say that $E'$ is a \emph{major subset} of $E$ if $E' \subset E$ and
$|E'| \geq \frac{1}{2}|E|$.
\end{definition}

\begin{definition} If $E$ is a set of finite measure, we let $X(E)$ denote the space of all functions $F$ supported on $E$ such that $\|F\|_\infty \leq 1$.
\end{definition}

\begin{definition}\label{rwt-def}  If $\alpha = (\alpha_1,\ldots, \alpha_n)$ 
is an admissible tuple, we say that an $n$-linear form $\Lambda$ is of 
restricted type $\alpha$ if for every 
sequence $E_1, \ldots, E_n$ of subsets of $\R$ with finite measure, 
there exists 
a major subset 
$E'_j$ of $E_j$ for each bad index $j$ (one or none) such that
\be{cut} | \Lambda(F_1, \ldots, F_n) | \lesssim
|E|^{\alpha}
\end{equation}
for all functions $F_i \in X(E'_i)$, $i=1, \ldots n$,
where we adopt the convention $E'_i = E_i$ for good indices $i$, 
and $|E|^{\alpha}$ is shorthand for
$$ |E|^{\alpha} = |E_1|^{\alpha_1} \ldots |E_n|^{\alpha_n}.$$
\end{definition}

The restricted type result we will prove directly is:
\begin{theorem}\label{rwt-bound}  
The form $\Lambda$ as in Theorem \ref{main} is of restricted 
type $\alpha$ for all bad tuples $\alpha$ such that 
\be{p-bound}
1/2 < \alpha_i < 1
\end{equation}
for all all good indices $i$ and
\be{lick}
k - n + 2 > \alpha_j > k - n + \frac{3}{2}
\end{equation}
for the bad index $j$.
\end{theorem}

Once Theorem \ref{rwt-bound} is granted, which we shall assume throughout
this section, the issue of proving 
Theorem \ref{main} is to pass to the
(admissible part of the) convex hull of tuples described in 
Theorem \ref{rwt-bound}
and convert restricted type estimates to strong type estimates.

We need the following easy lemma on permutahedrons:
\begin{lemma}  \label{honey}
Let $a_1 > \ldots > a_n$ be numbers.  Then the convex hull of all permutations of $(a_1, \ldots, a_n)$ consists of those points $(x_1, \ldots, x_n)$ such that $x_1 + \ldots + x_n = a_1 + \ldots + a_n$ and
$x_{i_1} + \ldots + x_{i_r} \leq a_1 + \ldots + a_r$
for all $1 \leq i_1 < \ldots < i_r \leq n$ and $1 \leq r \leq n$.
\end{lemma}

\begin{proof}  It is clear that the convex hull belongs to the set described above.  It thus suffices to show that the only extreme points of the above set are the permutations of $(a_1, \ldots, a_n)$.

Let $(x_1, \ldots, x_n)$ be an extreme point; by symmetry we may assume that
$x_1 \geq \ldots \geq x_n$.  If $x_1 + \ldots + x_r < a_1 + \ldots + a_r$ for
some $1 \leq r < n$, then we may modify $x$ by a small multiple of $e_r - e_{r+1}$ in either direction without leaving the set, contradicting the extremality of $x$.  Thus
$x_1 + \ldots + x_r = a_1 + \ldots + a_r$ for all $r$, so that
$(x_1, \ldots, x_n) = (a_1, \ldots, a_n)$, as desired.
\end{proof}

Now let $P$ denote the set of all admissible tuples described by
Theorem \ref{rwt-bound} and let $Q$ denote the set of all admissible
tuples $\alpha$ such that
\begin{equation}\label{region-q}
\alpha_{i_1}+\dots+\alpha_{i_r}<\frac {n-2k+r}{2}
\end{equation}
for all $1\le i_1<\dots <i_r\le n$ and $1\le r\le n$.

We have
\begin{lemma}  \label{convexhull}
The set $Q$ is contained in the convex hull of $P$.
The set $Q$ contains all good tuples.
If $\alpha\in Q$ has bad index $j$, then there is a
$\tilde{\alpha} \in P$ with $\tilde{\alpha}_i>\alpha_i$ for all
$i\neq j$ and such that $\alpha$ is in the convex hull of $\tilde{\alpha}$
and the elements in $Q$ whose bad index is not equal to $j$.

\end{lemma}

\begin{proof}  
For the first statement it suffices to prove that all tuples satisfying
\eqref{admissible}, \eqref{admissum}, and \eqref{region-q} are
contained in the convex hull of $\overline{P}$. This in turn follows
immediately from Lemma \ref{honey} and the observation that $\overline{P}$
contains all tuples which have $n-2k$ elements equal to $1$,
$2k-1$ elements equal to $1/2$, and the remaining element equal to
$3/2+k-n$.

The second statement follows immediately from the observation that the 
right hand side of (\ref{region-q}) is greater than or equal to one with strict
inequality in case $r>1$.

To see the third statement, assume by symmetry that
$\alpha\in Q$ has bad index $n$.
Define $\tilde{\alpha}_i=\max({\alpha_i},1/2)$ if $i\neq n$
and $\tilde{\alpha}_n=3/2+k-n$. Then \eqref{region-q} shows
that $\sum_i\tilde{\alpha_i}<1$, so we can enlarge the entries
of $\tilde{\alpha}$ so that $\tilde{\alpha}\in P$ and
$\tilde{\alpha}_i>\alpha_i$ for $i\neq n$.
We can write $\alpha=\theta \tilde{\alpha}+(1-\theta)\alpha'$
where $\alpha'$ is some tuple with $\alpha'_n=1/2$. 
Then a similar application of
Lemma \ref{honey} as before implies that
$\alpha'$ is in the convex hull of those elements $\alpha''\in \overline{P}$
for which $\alpha_n''=1/2$. This implies the third statement of Lemma 
\ref{convexhull}.

\end{proof}

We first discuss good exponent $n$-tuples:

\begin{lemma}\label{goodtuples}  
Let the assumptions and notation be as in Theorem
\ref{main}.
Then $\Lambda$ is of restricted type $\alpha$ for all good 
tuples $\alpha$. 
\end{lemma}

\begin{proof}
Let $\alpha$ be a good tuple. By symmetry we can assume 
that $\alpha_n=\max_i \alpha_i$. Then we have
$\alpha_n>0$ and $\alpha_i\le 1/2$ for all $i\neq n$.
By Lemma \ref{convexhull} we find $\theta_j\ge 0$ such that
\begin{equation}\label{conv}
{\alpha}=\sum_{j=1}^n \theta_j \alpha^{(j)},\ \ \ \sum_{i=1}^n\theta_j=1
\end{equation}
where $\alpha^{(j)}=(\alpha^{(j)}_i)_{i=1}^n$ is an admissible tuple in 
$P$ with bad index $j$. We can arrange that $\theta_n>0$.

For $\lambda>0$ let $A(\lambda)$ be the best constant such that
$$
 |\Lambda(F_{1}, \ldots, F_n)| \le  A(\lambda)  |E|^{\alpha}\ \ \ .
$$
for all 
sets $E_1, \ldots, E_n$ of finite measure with
\begin{equation}\label{lambda}
|E|^{\alpha^{(n)}}<\lambda |E|^{\alpha} 
\end{equation}
and functions 
$F_i\in X(E_i)$.
Let $A(\infty)$ be the supremum of all $A(\lambda)$.
By the a priori smoothness and support assumptions on 
$m$, $A(\infty)$ is finite
and the point is to prove that it is bounded.

By splitting
$$
 |\Lambda(F_{1}, \ldots, F_n)| \le
 |\Lambda(F_{1}, \ldots, F_n\chi_{E'_n})| 
+ |\Lambda(F_{1}, \ldots, F_n\chi_{E_n\setminus E'_n})| 
$$
appropriately and using restricted type $\alpha^{(n)}$ from Theorem \ref{rwt-bound}
we obtain
\begin{equation}\label{firstalambda}
A(\lambda)\le C \lambda + 2^{-\alpha_n}A(\infty)\ \ \ .
\end{equation}
On the other hand, if 
$$
\frac \lambda 2 |E|^{\alpha} \le |E|^{\alpha^{(n)}}<\lambda |E|^{\alpha} \ \ \ ,
$$
then we can use (\ref{conv}) to find an index $j\neq n$ such that
$$
|E|^{\alpha^{(j)}}\le C\lambda^{-\frac {\theta_n}{1-\theta_n}}|E|^{\alpha}\ \ \ .
$$
By splitting
$$
 |\Lambda(F_{1}, \ldots, F_n)| \le
 |\Lambda(F_{1}, \ldots, F_j\chi_{E'_j},\dots,F_n)| 
+ |\Lambda(F_{1}, \ldots, F_j\chi_{E_j\setminus E'_j},\dots,F_n)| 
$$
appropriately and using restricted type $\alpha^{(j)}$ from Theorem \ref{rwt-bound}
we obtain
$$
A(\lambda)\le \max(A(\lambda/2),C\lambda^{-\frac {\theta_n}{1-\theta_n}}+A(\max_{j\neq n} 
2^{\alpha_j-\alpha_j^{(n)}}\lambda))\ \ \ .
$$
Since $\max_{j\neq n}(\alpha_j-\alpha_j^{(n)})$ is negative, we can iterate
the previous inequality to obtain for sufficiently large $\lambda$:
$$
A(\infty)\le 1+ A(\lambda)\ \ \ .
$$
Combining this with (\ref{firstalambda}) gives
$$
A(\infty)\le C+ 2^{-\alpha_n}A(\infty)\ \ \ ,
$$
which proves boundedness of $A(\infty)$.
\end{proof}

\begin{lemma}\label{good-tuples}
Let $1<p_i\le \infty$ for $1\le i\le n$ such that \eqref{scaling}
holds. Then
$$
\Lambda(f_1,\dots,f_n)\le C \|f_1\|_{p_1}\dots \|f_n\|_{p_n}
$$
for all functions $f_i$ supported on a set of finite measure.
\end{lemma}

\begin{proof}
By symmetry we can assume that $p_i\neq \infty$ for $i\le j$
and $p_i=\infty$ for $i>j$ for a certain $j$.
Lemma \ref{goodtuples} implies
$$
\Lambda(f_1,\dots f_n)\le C\|f_1\|_{L^{q_1,1}}\dots 
\|f_j\|_{L^{q_j,1}} \|f_{j+1}\|_\infty
\|f_{n}\|_\infty
$$
for all $q_1,\dots q_j$ in a small neighborhood of $p_1,\dots ,p_j$
satisfying
$$
1/q_1+\dots 1/q_j=1\ \ \ .
$$
Fix functions $f_{j+1},\dots , f_n$. Then
Marcinkiewicz interpolation as in \cite{janson} implies
$$
\Lambda(f_1,\dots f_n)\le C\|f_1\|_{L^{p_1}}\dots 
\|f_j\|_{L^{p_j}} \|f_{j+1}\|_\infty
\|f_{n}\|_\infty
$$
for all functions $f_1$, \dots $f_j$.
\end{proof}

We turn to bad tuples $\alpha$.

\begin{lemma}\label{badtuples} 
Let the assumptions and notation be as in Theorem
\ref{main}. Then $\Lambda$ is of restricted type $\alpha$ for all bad  
tuples $\alpha$ satisfying (\ref{region-q}).
\end{lemma}

\begin{proof}  Fix $\alpha = (\alpha_1, \ldots, \alpha_n)$; by symmetry 
we may assume that $\alpha$ has bad index $n$.
By Lemma \ref{convexhull} we find $\theta_j$ such that
\begin{equation}\label{conv2}
{\alpha}=\sum_{j=1}^n \theta_j \alpha^{(j)},\ \ \ \sum_{i=1}^n\theta_j=1
\end{equation}
where $\alpha^{(j)}=(\alpha^{(j)}_i)_{i=1}^n$ is an admissible tuple in 
$P$ with bad index $j$. We have $\theta_n>0$ and by the last statement of
Lemma \ref{convexhull}  we can assume that
$\alpha_j-\alpha_j^{(n)}$ is negative for all $j\neq n$.

For $\lambda>0$ let $A(\lambda)$ be the best constant such that
for all 
sets $E_1, \ldots, E_n$ of finite measure with
\begin{equation}\label{lambda2}
|E|^{\alpha^{(n)}}<\lambda |E|^{\alpha} 
\end{equation}
there is a major subset $E_n'$ of $E_n$ such that
$$
 |\Lambda(F_{1}, \ldots, F_n)| \le  A(\lambda)  |E|^{\alpha}\ \ \ .
$$
for all 
$F_i\in X(E_i')$.
Let $A(\infty)$ be the supremum of all $A(\lambda)$.

Using restricted type $\alpha^{(n)}$ from Theorem \ref{rwt-bound} we obtain
\begin{equation}\label{alcl}
A(\lambda)\le C\lambda
\end{equation}

On the other hand let
$$
\frac \lambda 2 |E|^{\alpha} \le |E|^{\alpha^{(n)}}<\lambda |E|^{\alpha} \ \ \ ,
$$
then we can find an index $j\neq n$ such that
$$
|E|^{\alpha_j}\le C\lambda^{-\frac {\theta_n}{1-\theta_n}}|E|^{\alpha}\ \ \ ,
$$
and we can use restricted type $\alpha^{(j)}$ from Theorem \ref{rwt-bound}
to conclude
$$
A(\lambda)\le \max(A(\lambda/2),C\lambda^{-\frac {\theta_n}{1-\theta_n}}+
A(\max_{j\neq n } 
2^{\alpha_j-\alpha_j^{(n)}}\lambda))\ \ \ .
$$
Since $\max_{j\neq n}(\alpha_j-\alpha_j^{(n)})$ is negative, we can iterate
the previous inequality to obtain for sufficiently large $\lambda$:
$$
A(\infty)\le 1+ A(\lambda)\ \ \ .
$$
Together with (\ref{alcl}) this proves bondedness of $A(\infty)$.

\end{proof}

Finally, we convert restricted type estimates for bad tuples
$\alpha$ into strong type estimates by proving a 
Marcinkiewicz interpolation result
in the spirit of \cite{janson}.

\begin{lemma}\label{bad-tuples}
Let $\alpha$ be a bad 
tuple satisfying (\ref{region-q}) and assume that $n$ is the bad index.
Set $p_i=1/\alpha_i$ for $1\le i \le n$.
Then
$$
\left\|T(f_1,\dots,f_{n-1})\right\|_{p_n'}
\le C \|f_1\|_{p_1}\dots \|f_{n-1}\|_{p_{n-1}}
$$
for all functions $f_i$ supported on a set of finite measure.
\end{lemma}

\begin{proof}

We assume for simplicity that $p_i\neq \infty$ for all $i$.
If this was not the case, we could freeze the function  $f_i$ 
and the exponent $p_i$ whenever $p_i=\infty$ and run the argument 
on the remaining functions only, as done in the proof of 
Lemma \ref{good-tuples}.

Let $f_1, \ldots, f_{n-1}$ be functions such that $\|f_i\|_{p_i} = 1$ for $1 \leq i < n$.  We have to show that
$$ \| T(f_1, \ldots, f_{n-1})\|_{p'_n} \lesssim 1.$$
We may assume that the $f_i$ are non-negative.  By a measure-preserving rearrangement, we may assume that the $f_i$ are supported on the half-line
$(0,\infty)$ and are monotone non-increasing on this half-line.

Let $\chi_k$ denote the function $\chi_k = \chi_{(2^k,2^{k+1}]}$.  We can expand
the desired estimate as
$$ \| \sum_{k_1, \ldots, k_{n-1}} 
T(f_1 \chi_{k_1}, \ldots, f_{n-1} \chi_{k_{n-1}}) \|_{p'_n} \lesssim 1.$$
Since $p'_n \leq 1$, we have the elementary inequality
$$ \| \sum_\beta F_\beta \|_{p'_n}^{p'_n} \leq \sum_\beta \| F_\beta\|_{p'_n}^{p'_n},$$
so it suffices to show that
\be{oak}
 \sum_{k_1, \ldots, k_{n-1}} \| 
T(f_1 \chi_{k_1}, \ldots, f_{n-1} \chi_{k_{n-1}})
\|_{p'_n}^{p'_n} \lesssim 1.
\end{equation}
By symmetry we may restrict the summation to the region
$$ k_1 \geq k_2 \geq \ldots \geq k_{n-1}.$$
Fix $k_1, \ldots, k_{n-1}$.  Let $\lambda > 0$ be arbitrary, and consider the set
$$ E_n = \{ \Re T(f_1 \chi_{k_1}, \ldots, f_{n-1} \chi_{k_{n-1}}) > \lambda \}.$$
Let $\alpha$ be an admissible tuple close to $1/p$; we may thus 
assume $\alpha$ has bad index $n$.  Since $\Lambda$ is of restricted type 
$\alpha$, and $f_i \chi_{k_i} \in f_i(2^{k_i}) X((2^k,2^{k+1}])$, we may thus 
find a major subset $E'_n$ of $E_n$ such that
$$ |\Lambda(f_1 \chi_{k_1}, \ldots, f_{n-1} \chi_{k_{n-1}}, \chi_{E'_n})|
\lesssim |E_n|^{\alpha_n} \prod_{i=1}^{n-1} f_i(2^{k_i}) 2^{k_i\alpha_i}.$$
By definition of $E_n$, we thus have
$$ \lambda |E_n| \lesssim |E_n|^{\alpha_n} \prod_{i=1}^{n-1} f_i(2^{k_i}) 
2^{k_i\alpha_i}.$$
Solving for $|E_n|$, and optimizing in $\alpha$, one obtains
$$ |E_n| \lesssim \lambda^{-p'_n}
2^{-\eps(k_1-k_{n-1})} \min(\frac{F}{\lambda}, \frac{\lambda}{F})^\eps
(\prod_{i=1}^{n-1} f_i(2^{k_i}) 2^{k_i/p_i})^{p'_n}
$$
for some $\eps > 0$, where $F = \prod_{i=1}^{n-1} f_i(2^{k_i})$.  By symmetry one may obtain the same bound when $E_n$ is replaced by
$$ \{ |T(f_1 \chi_{k_1}, \ldots, f_{n-1} \chi_{k_{n-1}})| > \lambda \}.$$
Integrating this over all $\lambda$, one then obtains
$$ 
\left\|T(f_1 \chi_{k_1}, \ldots, f_{n-1} \chi_{k_{n-1}})
\right\|_{p'_n} \lesssim 2^{-\eps(k_1 - k_{n-1})} \prod_{i=1}^{n-1} f_i(2^{k_i}) 2^{k_i/p_i}.
$$
To prove \eqref{oak}, it thus suffices to show
\be{oak2} (\sum_{k_1 \geq \ldots \geq k_{n-1}}
2^{-\eps(k_1 - k_{n-1})} (\prod_{i=1}^{n-1} f_i(2^{k_i}) 2^{k_i/p_i})^{p'_n})^{1/p'_n} \lesssim 1.
\end{equation}
Write $s = k_1 - k_{n-1}$.  For fixed $s$ and $k_1$ there are at most $(1+s)^C$
choices of $k_i$.  Fixing $s$, and then applying H\"older's inequality using
\eqref{scaling}, we can
estimate the left-hand side of \eqref{oak2} by
$$ \sum_{s \geq 0} (1+s)^C 2^{-\eps s}
\prod_{i=1}^{n-1} (\sum_k (f_i(2^k) 2^{k/p_i})^{p_i})^{1/p_i}.$$
The $s$ sum is convergent, and the expression inside the product is essentially
$\|f_i\|_{p_i} = 1$.  The claim is thus proved.
\end{proof}

Theorem \ref{main} now follows from
Lemma \ref{good-tuples} and Lemma \ref{bad-tuples}.

\section{Exceptional set}\label{except}

It remains to prove Theorem \ref{rwt-bound}.  Let $p$ satisfy the hypotheses of the theorem; by symmetry we may assume that the bad index of $p$ is $n$.  We have to show that for any $E_1, \ldots, E_n$ one can find a major subset $E'_n$ of $E_n$ such that \eqref{cut} holds for all $F_i \in X(E'_i)$.  By \eqref{scaling} and a scaling argument one may take $|E_n| = 1$.  

We shall define $E'_n$ explicitly as
\be{en-def} E'_n = \{ x \in E_n: \M \chi_{E_i}(x) < C |E_i| \hbox{ for all } 1 \leq i \leq n \}.
\end{equation}
From the Hardy-Littlewood maximal inequality we see that $|E_n \backslash E'_n| \leq \frac{1}{2}$ if
$C$ is chosen sufficiently large.  Thus we have $|E'_n| \geq \frac{1}{2} |E_n|$
as desired.

Let $F_i$ be arbitrary elements of $X(E'_i)$.
Define the normalized functions $f_1, \ldots, f_n$ by
$$ f_i = \frac{F_i \chi_{E'_i}}{|E'_i|^{1/2}}, \quad i = 1, \ldots n;$$
note that 
\be{l2-norm}
\| f_i \|_2 \lesssim 1 \hbox{ for all } i = 1, \ldots, n.
\end{equation}
Also define the numbers
$a_i = |E_i|^{1/2}$.  We may rewrite \eqref{cut} as
\be{lam}
 |\Lambda(f_1, \ldots, f_n)| \lesssim \prod_{i=1}^n a_i^{\theta_i}
\end{equation}
where $\theta_i = \frac{2}{p_i}-1$ for $1 \leq i \leq n-1$.  Since $a_n = 1$, the value of $\theta_n$ is arbitrary, but we shall set it so that
\be{theta}
 \theta_1 + \ldots + \theta_n = n-2k.
\end{equation}
From \eqref{scaling}, \eqref{p-bound} and \eqref{lick} we see that
\be{theta-bound}
 0 < \theta_i < 1 \hbox{ for all } i = 1, \ldots, n.
\end{equation}

Let $N \gg 1$ be a large constant to be chosen later.
For any interval $I$ and $1 \leq i \leq n$, define the normalized averages $\lambda_i(I)$ by
$$ \lambda_i(I) = \frac{1}{|I| |E'_i|} \int_{E'_i} \tilde \chi_I^N.$$
Clearly we have the estimates
\be{f1}
\| f_i \tilde \chi_I^N \|_1 \lesssim a_i \lambda_i(I) |I|
\end{equation}
and
\be{f2}
\| f_i \tilde \chi_I^{N/2} \|_2 \lesssim \lambda_i(I)^{1/2} |I|^{1/2}
\end{equation}
for all $I$ and $i$.  

From the construction of $E'_n$ we see that the $\lambda_i$ cannot simultaneously be large.  More precisely, we have

\begin{lemma}  For any interval $I$ we have
\be{decay}
\lambda_n(I) \lesssim (1 + \lambda_1(I) + \ldots + \lambda_{n-1}(I))^{1-N}.
\end{equation}
\end{lemma}

\begin{proof}
Suppose first that $2I$ intersected $E'_i$.  Then there exists $x \in 2I$
such that $\M \chi_{E_i}(x) \leq |E_i|$ for all $1 \leq i \leq n$.  This implies
that $\lambda_i(I) \lesssim 1$ for all $1 \leq i \leq n$, which implies \eqref{decay}.

Now suppose that $j\geq 1$ was such that $2^j I$ was disjoint from $E'_i$,
but $2^{j+1} I$ intersected $E'_i$.  By arguing as before we see that
$\lambda_i(I) \lesssim 2^j$ for all $1 \leq i \leq n-1$, and
$\lambda_n(I) \lesssim 2^{-j(1-N)}$, which again implies \eqref{decay}.
\end{proof}

As we shall see, the dominant contribution to \eqref{lam} shall come from those intervals $I$ for which $\lambda_i(I) \sim 1$.

To prove Theorem \ref{rwt-bound} it thus suffices to prove the following estimate.

\begin{theorem}\label{abstract}  Let $\Lambda$ be as above, let $f_1, \ldots, f_n$ be functions satisfying \eqref{l2-norm}, and $a_1, \ldots, a_n$ be positive numbers.
For each interval $I$ and $1 \leq i \leq n$ we let $\lambda_i(I)$ be a non-negative number such that \eqref{f1}, \eqref{f2}, \eqref{decay} hold for all $I$ and $i$.  Then for any $\theta_i$ satisfying \eqref{theta} and \eqref{theta-bound} we have \eqref{lam}, provided that $N$ is chosen sufficiently large depending on $\theta$.
\end{theorem}

We have thus reduced the problem to that of estimating $\Lambda$ on functions which are $L^2$-normalized, and whose $L^1$ and $L^2$ averages on intervals are somewhat under control.

\section{Discretization}\label{discrete}

Let $f_i$, $a_i$, $\lambda_i(I)$ be as in Theorem \ref{abstract}.
We now decompose the multiplier $m$ using a Whitney decomposition, and replace $\Lambda$ with a discretized variant.

We may extend $m$ from the $n-1$-dimensional hyperplane $\hyperplane$ to the
entire space $\R^n$ in such a way that \eqref{symb-gamma} holds for all $\xi \in \R^n \backslash \subspace$ and all derivatives $\alpha$ up to a sufficiently large order.

Define a \emph{shifted $n$-dyadic mesh} $D = D^n_{\alpha}$ to be a collection of cubes of the form
$$ D^n_{\alpha} = \{ 2^j (k + (0,1)^n + (-1)^j \alpha):
: j \in \Z, \quad k \in \Z^n \}$$
where $\alpha \in \{ 0, \frac{1}{3}, \frac{2}{3} \}^n$.  We define a \emph{shifted dyadic cube} to be any member of a shifted $n$-dyadic mesh.

Observe that for every cube $Q$, there exists a shifted dyadic cube\footnote{This observation is due to Michael Christ.} $Q'$ such that $Q \subseteq \frac{9}{10} Q'$ and $|Q'| \sim |Q|$; this is best seen by first
verifying the $n=1$ case.

Consider the collection $\Q$ of all shifted dyadic cubes $Q$ such that
$$ \dist(Q,\subspace) \sim C_0 \diam(Q);$$
here $C_0$ is a large constant to be chosen later.
From the above observation we see that the cubes $\{\frac{9}{10} Q: Q \in \Q\}$ form a finitely overlapping cover of $\R^n \backslash \subspace$.  
This implies that we may partition
\be{m-decom} m = \sum_{Q \in \Q} m_Q
\end{equation}
where each $m_Q$ is supported in $Q\cap \hyperplane$ and satisfies the bounds
\be{symb-local}
| \partial_\xi^\alpha m_Q(\xi)| \lesssim \diam(Q)^{-|\alpha|}
\end{equation}
for all derivatives $\partial_\xi^\alpha$ on $\Gamma$ up to some 
sufficiently large order.

From \eqref{m-decom} we have
$$
\Lambda = \sum_{Q \in \Q} \Lambda_{m_Q}.
$$
Of course $\Lambda_{m_Q}$ vanishes unless $Q$ intersects $\hyperplane$.
Since there are only a finite number of shifted dyadic meshes, we see that \eqref{lam} will follow from
$$
\sum_{Q \in \Q \cap D, Q \cap \hyperplane \neq \emptyset}
|\Lambda_{m_Q}(f_1, \ldots, f_n)| \lesssim \prod_{i=1}^n a_i^{\theta_i}
$$
where $D = D^n_\alpha$ is any shifted dyadic mesh.  Henceforth
$\alpha = (\alpha_1, \ldots, \alpha_n)$ will be fixed.

To estimate the contribution of each $\Lambda_Q$ we introduce tiles in
the time-frequency plane $\R \times \R$.

\begin{definition}  Let $1 \leq i \leq n$.  An $i$-tile is a rectangle
$P = I_P \times w_P$ with area 1 and with $I_P \in D^1_0$, $w_P \in D^1_{\alpha_i}$.  A multi-tile is an $n$-tuple $\pv = (P_1, \ldots, P_n)$
such that each $P_i$ is an $i$-tile, and the $I_{P_i} = I_\pv$ are independent 
of $i$.  The frequency cube $Q_\pv$ of a multi-tile is defined to be
$\prod_{i=1}^n w_{P_i}$.
\end{definition}

If $\pv$ appears in an expression, we shall always adopt the convention that $P_i$ denotes the $i^{th}$ component of $\pv$.

\begin{definition}  Let $1 \leq i \leq n$, and let $P$ be an $i$-tile.  The semi-norm $\| f\|_P$ is defined by
$$ \|f\|_P = \frac{1}{|I_P|} \| (\Delta_{w_P} f) \tilde \chi_{I_P}^{2N} \|_1,$$
where $\Delta_{w_P}$ is a Fourier multiplier whose symbol $\psi_{w_P}$ is a bump function adapted to $w_P$ and which equals 1 on $\frac{9}{10} w_P$.  
\end{definition}

The quantity $\|f\|_P$ can be viewed as an average value of $f$ on the time-frequency tile $P$.  From the rapid decay of $\Delta_{w_P} f$ we observe the
crude estimate

\begin{lemma}\label{triv}  For any $P$, we have
$$ \| f\|_P \lesssim \frac{1}{|I_P|} \| f \tilde \chi_{I_P}^{2N} \|_1.$$
\end{lemma}

The relationship between these semi-norms and the $\Lambda_{m_Q}$ is given by

\begin{lemma}\label{spatial} For any $Q \in \Q \cap D$, we have
$$ |\Lambda_{m_Q}(f_1, \ldots, f_n)|
\lesssim \sum_{\pv: Q_\pv = Q} |I_\pv| \prod_{i=1}^n \| f_i \|_{P_i}$$
where $\pv$ runs over all multi-tiles with frequency cube $Q$.
\end{lemma}

\begin{proof}  By translation and scale invariance we may make $Q$ the unit cube $[0,1]^n$.

We may write $m_Q(\xi) = \tilde m(\xi) \prod_{i=1}^n \psi_{w_{P_i}}(\xi_i)$, where $\tilde m$ is supported on $[0,1]^n\cap \hyperplane$ 
and satisfies the same bounds \eqref{symb-local} as $m_Q$; in other words, $\tilde m$ is a bump function
on $\hyperplane$.  Since
$$ \Lambda_{m_Q}(f_1, \ldots, f_n) = \Lambda_{\tilde m}(\Delta_{w_{P_1}} f_1,
\ldots \Delta_{w_{P_n}} f_n),$$
it suffices to show the estimate
$$ |\Lambda_{\tilde m}(g_1, \ldots, g_n)| \lesssim
\sum_l \prod_{i=1}^n \| g_i \tilde 
\chi_{[l,l+1]}^N \|_1.$$
From Plancherel's theorem and \eqref{symb-local} one sees that
$$ \Lambda_{\tilde m}(g_1, \ldots, g_n) = \int K(x) \prod_{i=1}^n g_i(x_i)\ dx$$
where $x = (x_1, \ldots, x_n)$ and the kernel $K$ satisfies the estimate
$$ |K(x)| \lesssim (1 + \sum_{i,j} |x_i - x_j|)^{-M}$$
for arbitrarily large $M$.  In particular, we have
$$ |K(x)| \lesssim \sum_l \prod_{i=1}^n \tilde \chi_{[l,l+1]}^{2N}(x_i)$$
and the claim follows.
\end{proof}

Let $\Pv$ denote the set of all multi-tiles $\pv$ such that $Q_{\pv} \in \Q \cap D$ and $Q_{\pv}$ intersects $\hyperplane$.  From the above lemma, it suffices to show that
\be{disc-wosh}
\sum_{\pv \in \Pv} |I_\pv| \prod_{i=1}^n \|f_i\|_{P_i}
\lesssim \prod_{i=1}^n a_i^{\theta_i}.
\end{equation}
Note that the multiplier $m$ no longer plays a role.

\section{Rank}\label{rank-sec}

The tiles in $\Pv$ have essentially $k$ independent frequency parameters.  To
make this more precise we need some notation.

\begin{definition}  Let $P$ and $P'$ be tiles.  We write $P' < P$ if $I_{P'} \subsetneq I_P$ and $w_P \subseteq 3w_{P'}$, and $P' \leq P$ if $P' < P$ or $P' = P$.
We write $P' \lesssim P$ if $I_{P'} \subseteq I_P$ and $w_P \subseteq C C_0 w_{P'}$.  We write $P' \lesssim' P$ if $P' \lesssim P$ and
$P' \not \leq P$.
\end{definition}

Note that the ordering $<$ is slightly different from the one in Fefferman \cite{fefferman}
or Lacey and Thiele \cite{laceyt1}, \cite{laceyt2}, \cite{thiele} as $P'$ and $P$ do not quite have to intersect.  This slightly less strict ordering is more
convenient for technical purposes.

If $C_0$ is sufficiently large, then we have

\begin{lemma}\label{rank}  Let $1 \leq i_1 < i_2 < \ldots < i_k \leq n$ be integers, and $\pv$, $\pv'$ be multi-tiles in $\Pv$.  If $P'_{i_s} \leq P_{i_s}$ for all $s = 1, \ldots, k$, then $P'_i \lesssim P_i$ for all $1 \leq i \leq n$.  If we further
assume that $|I_{\pv'}| \ll |I_{\pv}|$, then we have $P'_i \lesssim' P_i$ for at least two choices of $i$.
\end{lemma}

\begin{proof}
Since $\subspace$ is non-degenerate, we can write it as a graph
$$ \{ \xi: \xi = h(\xi_{i_1}, \ldots, \xi_{i_k}) \},$$
where $h$ is a linear map from $\R^k$ to $\hyperplane$.

Let $\xi$, $\xi'$ denote the centers of $Q_\pv$ and $Q_{\pv'}$ respectively. From the definition of $\Pv$ we have
\be{xid}
|\xi - h(\xi_{i_1}, \ldots, \xi_{i_k})| \sim C_0 |I_\pv|^{-1}.
\end{equation}
and
\be{shum} 
|\xi_1 + \ldots + \xi_n| \lesssim |I_\pv|^{-1},
\end{equation}
and similarly for $\xi'$.  Since $3w_{P'_{i_s}}$ contains $w_{P_{i_s}}$, we have
$$ \xi_{i_s} = \xi'_{i_s} + O(|I_{\pv'}|^{-1}).$$
Combining this with \eqref{xid} we see that
$$ \xi = \xi' + O(C_0 |I_{\pv'}|^{-1}),$$
which implies that $P'_i \lesssim P_i$ for all $I$ as desired.

Now suppose $|I_{\pv'}| \ll |I_{\pv}|$.  By subtracting \eqref{xid} for $\xi$ and $\xi'$ we thus have
$$
|(\xi-\xi') - h((\xi - \xi')_{i_1}, \ldots, (\xi - \xi')_{i_k})| \sim C_0 |I_{\pv'}|^{-1},$$
which implies that
$$ |\xi - \xi'| \gtrsim C_0 |I_{\pv'}|^{-1}.$$
On the other hand, from \eqref{shum} we have
$$\left| (\xi - \xi')_1 + \ldots (\xi - \xi')_n \right| \sim |I_{\pv'}|^{-1}.$$
If $C_0$ is sufficiently large, this guarantees that there exist $1 \leq i < i' \leq n$ such that
$$ |(\xi - \xi')_i|, |(\xi - \xi')_{i'}| \geq 3 |I_{\pv'}|^{-1},$$
which combined with the previous observations gives $P'_i \lesssim' P_i$ and
$P'_{i'} \lesssim' P_{i'}$ as desired.
\end{proof} 

\begin{definition}  If $\Pv$ is a collection of tiles, we define the norm
$\| f_i \|_{\Pv,i}$ by
$$\| f_i \|_{\Pv,i} = \sup_{\pv \in \Pv} \| f_i \|_{P_i}$$
\end{definition}

We now claim that Theorem \ref{abstract} follows from 

\begin{theorem}\label{discretized} Let $f_1, \ldots, f_n$ be functions obeying \eqref{l2-norm}, and $a_1, \ldots, a_n$, $\lambda_1, \ldots, \lambda_n$ be positive numbers.  Let $\Pv$ be a finite collection of multi-tiles such that Lemma \ref{rank} holds, and such that
\be{F1}
\| f_i \tilde \chi_{I_\pv}^N \|_1 \lesssim a_i \lambda_i |I_\pv|,
\end{equation}
\be{F2}
\| f_i \tilde \chi_{I_\pv}^{N/2} \|_2 \lesssim \lambda_i^{1/2} |I_\pv|^{1/2},
\end{equation}
for all $\pv \in \Pv$ and $1 \leq i \leq n$.  Let $I_0$ be an interval such that
$I_\pv \subseteq I_0$ for all $\pv \in \Pv$, and
\be{i0}
\| f_i \tilde \chi_{I_0}^{N/2} \|_2 \lesssim \lambda_i^{1/2} |I_0|^{1/2},
\end{equation}
for all $1 \leq i \leq n$.  Then one has
\be{disc}
\sum_{\pv \in \Pv} |I_\pv| \prod_{i=1}^n \|f_i\|_{P_i}
\lesssim A^{n-2k-\theta_1 - \ldots - \theta_n} \min(1, |I_0|)
\prod_{i=1}^n (\lambda_i a_i)^{\theta_i} (1 + \lambda_i),
\end{equation}
for any $\theta_i$ satisfying \eqref{theta-bound} and 
\be{theta-ineq}
\theta_1 + \ldots + \theta_n \leq n-k,
\end{equation}
where $A$ is the quantity
\be{a-def}
A = \sup_{1 \leq i \leq n} \|f_i\|_{\Pv,i}.
\end{equation}
\end{theorem}

Theorem \ref{discretized} contains some rather technical assumptions which are convenient for induction purposes.  In applications, we would only use the following corollary:

\begin{corollary}  Let $f_1, \ldots, f_n$, $a_1, \ldots, a_n$, $\lambda_1, \ldots, \lambda_n$, and $\Pv$ be as in the previous Theorem.  Then
$$ \sum_{\pv \in \Pv} |I_\pv| \prod_{i=1}^n \|f_i\|_{P_i}
\lesssim 
\prod_{i=1}^n (\lambda_i a_i)^{\theta_i} (1 + \lambda_i),
$$
for any $\theta_i$ satisfying \eqref{theta} and \eqref{theta-bound}.
\end{corollary}

Now let $\lambda_1, \ldots, \lambda_n$ be dyadic numbers such that
$$ \lambda_n \lesssim (1 + \lambda_1 + \ldots + \lambda_{n-1})^{1-N},$$
and apply the Corollary to those multi-tiles $\pv$ such that
$\lambda_i(P_i) \sim \lambda_i$ for $1 \leq i \leq n$.  The estimate \eqref{disc-wosh} then follows by summing in $\lambda_n$ and then in each of the $\lambda_i$, $1 \leq i \leq n-1$.

It remains to prove Theorem \ref{discretized}.  This shall be done in two stages.  Firstly we shall handle the case $k=1$, by arguments similar to those in Lacey and Thiele \cite{laceyt1}, \cite{laceyt2}, \cite{thiele}; this is the 
longest part of the proof, occupying Sections \ref{trees}-\ref{conclusion}.  
Then, in Section \ref{induct-sec}, we induct on $k$ to obtain the general case.

\section{Trees}\label{trees}

Let $k=1$.  Fix the $f_i$, $a_i$, $\lambda_i$, $\Pv$, and $I_0$. 

In order to estimate \eqref{disc} we shall have to organize $\Pv$ into trees, as in \cite{fefferman}, \cite{laceyt1}, \cite{laceyt2}, \cite{thiele}.

\begin{definition} For any $1 \leq j \leq n$ and a multi-tile $\pv_T \in \Pv$, define a $j$-tree with top $\pv_T$ to be a collection of multi-tiles $T \subseteq \Pv$ such that
$$ P_j \leq P_{T,j} \hbox{ for all } \pv \in T,$$
where $P_{T,j}$ is the $j$ component of $\pv_T$.  We write $I_T$ and $w_{T,j}$ for $I_{\pv_T}$ and $w_{P_{T,j}}$ respectively.  We say that $T$ is a tree if it is a $j$-tree for some $1 \leq j \leq n$.
\end{definition}

Note that $T$ does not necessarily have to contain its top $\pv_T$.

\begin{definition} 
For any tree $T$, define the $i$-size $\size_i(T)$ of $T$ to be the quantity
\be{size-def}
\size_i(T) = 
\left(\frac{1}{|I_T|} \sum_{\pv \in T: P_i \lesssim' P_{T,i}} |I_\pv| \| f_i\|_{P_i}^2\right)^{1/2} + \| f_i \|_{T,i}
\end{equation}
\end{definition}

The relationship between the $i$-size to \eqref{disc} is given by

\begin{lemma}\label{Cauchy-Schwarz}
If $T$ is a tree, then
\be{tree-bound}
\sum_{\pv \in T} |I_\pv| \prod_{i=1}^n \|f_i\|_{P_i}
\lesssim |I_T| \sup_{1 \leq i_1 < i_2 \leq n}
\size_{i_1}(T) \size_{i_2}(T) \prod_{i \neq i_1,i_2} \| f_i \|_{T,i}.
\end{equation}
\end{lemma}

\begin{proof}  We first deal with the contribution of those multi-tiles $\pv$ such that $|I_\pv| \sim |I_T|$.  From Lemma \ref{rank} there are only $O(1)$ of these multi-tiles, and the contribution can be handled by the estimate 
\be{sup}
\|f_i\|_{P_i} \leq \|f_i\|_{T,i}\leq \size_i(T).
\end{equation}
Now let us consider those multi-tiles for which $|I_\pv| \ll |I_T|$.  From Lemma \ref{rank} there exist $i_1$, $i_2$ such that $P_{i_s} \lesssim' P_{T,i_s}$
for $s=1,2$; by pigeonholing we may make $i_1$, $i_2$ independent of $\pv$.  If
one then uses \eqref{sup} for all $i \neq i_1, i_2$, one reduces to showing that
$$ \sum_{\pv \in T} |I_\pv| \| f_{i_1} \|_{P_{i_1}}
\| f_{i_2} \|_{P_{i_2}} \lesssim |I_T| \size_{i_1}(T) \size_{i_2}(T).
$$
But this follows from Cauchy-Schwarz.
\end{proof}

To apply Lemma \ref{Cauchy-Schwarz} we need to partition $\Pv$ into trees $T$ in such a way that we have good control on the $i$-sizes $\size_i(T)$ and the spatial sizes $|I_T|$.  This shall be done in four stages.

Firstly, in Section \ref{disjoint-sec}, we control the number of trees of a certain size by the following lemma.

\begin{definition} Let $1 \leq i \leq n$.  Two trees $T$, $T'$ are said to be \emph{strongly $i$-disjoint} if 
\bi
\item $P_i \neq P'_i$ for all $\pv \in T$, $\pv' \in T'$.
\item Whenever $\pv \in T$, $\pv' \in T'$ are such that
$w_{P_i} \subsetneq w_{P'_i}$, then one has
$I_{\pv'} \cap I_T = \emptyset$, and similarly with $T$ and $T'$ reversed.
\end{itemize}
\end{definition}

Note that if $T$ and $T'$ are strongly $i$-disjoint, then $P_i \cap P'_i = \emptyset$ for all $\pv \in T$, $\pv' \in T'$.

\begin{lemma}\label{disjoint-lemma}  Let $1 \leq i \leq n$, $m \in \Z$, and let $\T$ be a collection of trees in $\Pv$ which are mutually strongly $i$-disjoint and such that
\be{size-lower} \size_i(T) \sim 2^{-m} \hbox{ for all } T \in \T.
\end{equation}
Let $I_0$ be an interval such that $I_T \subseteq I_0$ for all $T \in \T$.
Then we have
\be{size-total} \sum_{T \in \T} |I_T| \lesssim 2^{2m} \| f_i \tilde \chi_{I_0}^{N/2} \|_2^2 \lesssim 2^{2m} \min(1, \lambda_i |I_0|).
\end{equation}
\end{lemma}

By applying Lemma \ref{disjoint-lemma} to singleton trees and $n=k=1$, one obtains

\begin{corollary}\label{lacey}  Let $f$ be a function, $m \in \Z$, $I_0$ be an interval, $\P$ be a collection of disjoint tiles such that $I_P \subseteq I_0$
and $\| f\|_P \sim 2^{-m}$ for all $P \in \P$.  Then we have
$$ \sum_{P \in \P} |I_P| \lesssim 2^{2m} \| f \tilde \chi_{I_0}^{N/2} \|_2^2.$$
\end{corollary}

In Section \ref{select-sec}, we use Lemma \ref{disjoint-lemma} to obtain the following tree selection algorithm.

\begin{lemma}\label{selection}  Let $1 \leq i \leq n$, $m \in \Z$, and suppose that one has
\be{m-bound} \size_i(T) \leq 2^{-m}
\end{equation}
for all trees $T$ in $\Pv$.  Then there exists a collection $\T$ of trees in $\Pv$ such that
\be{t-size}
 \sum_{T \in \T} |I_T| \lesssim 2^{2m} \min(1, \lambda_i |I_0|)
\end{equation}
and
\be{remainder} \size_i(T') \leq 2^{-m-1}
\end{equation}
for all trees $T'$ in $\Pv - \bigcup_{T \in \T} T$.
\end{lemma}

In Section \ref{lambda-sec}, we shall bound the $i$-size by

\begin{lemma}\label{lambda-lemma}  For any tree $T$ in $\Pv$ and $1 \leq i \leq n$, we have
$$ \size_i(T) \lesssim a_i \lambda_i.$$
\end{lemma}

Finally, in Section \ref{conclusion} we combine Lemma \ref{selection} and Lemma \ref{lambda-lemma} with Lemma \ref{Cauchy-Schwarz} to prove \eqref{disc} in the $k=1$ case.

\section{Proof of Lemma \ref{disjoint-lemma}}\label{disjoint-sec}

The second inequality in \eqref{size-total} follows from \eqref{F2} and the $L^2$-normalization of $f_i$, so it suffices to prove the first inequality.

Fix $1 \leq i \leq n$.  
By refining the trees $T$, we may assume that the tiles $\{ P_i: \pv \in T\}$ are all disjoint, and that
$$
\sum_{\pv \in T} |I_\pv| \| f_i \|_{P_i}^2 \sim 2^{-2m} |I_T|.
$$
In particular, we have
\be{albedo}
\sum_{\pv \in \bigcup_T T} |I_\pv| \| f_i \|_{P_i}^2 \sim 2^{-2m} \sum_T |I_T|.
\end{equation}
Also, from \eqref{size-lower} we have
\be{sup-size}
\| f_i \|_{P_i} \lesssim 2^{-m}
\end{equation}
for all $\pv \in \bigcup_T T$.

We shall shortly prove the estimate
\be{quasi-ortho}
\sum_{\pv \in \bigcup_T T} |I_\pv| \| f_i \|_{P_i}^2
\lesssim 
2^{-m} \| f_i \tilde \chi_{I_0}^{N/2} \|_2 (\sum_T |I_T|)^{1/2};
\end{equation}
the claim then follows by combining \eqref{albedo} and \eqref{quasi-ortho}.

The estimate \eqref{quasi-ortho} is somewhat reminiscent of an orthogonality estimate.  Accordingly, we shall use $TT^*$ methods and similar techniques in the proof.

By duality we may find a function $\phi_\pv$ for each $\pv \in \bigcup_T T$
such that $|\phi_\pv(x)| \lesssim \tilde \chi_{I_\pv}^{2N}(x)$ for all
$x \in \R$, and
$$ \| f_i \|_{P_i} = \frac{1}{|I_\pv|}
\langle \Delta_{w_{P_i}}^* \phi_\pv, f_i \rangle.$$
We can thus write the left-hand side of \eqref{quasi-ortho} as
$$ \langle \sum_{\pv \in \bigcup_T T} \| f_i \|_{P_i} \Delta_{w_{P_i}}^* \phi_\pv, f_i \rangle.$$
From the Cauchy-Schwarz inequality, the inequality \eqref{quasi-ortho} will follow from the estimate
\be{quasi-partial}
\| \sum_{\pv \in \bigcup_T T} \| f_i \|_{P_i} \Delta_{w_{P_i}}^* \phi_\pv
\tilde \chi_{I_0}^{-N/2} \|_2^2 \lesssim
2^{-2m} \sum_T |I_T|.
\end{equation}

Let us first consider the portion of the $L^2$ norm in \eqref{quasi-partial}
outside of $2I_0$.  From the triangle inequality, it will suffice to show that
\be{quasi-outside}
\| \sum_{\pv \in \bigcup_T T: I_\pv = I} 
\| f_i \|_{P_i} \Delta_{w_{P_i}}^* \phi_\pv \tilde \chi_{I_0}^{-N/2} \|_{L^2(\R \backslash 2I_0)}^2 \lesssim \frac{|I|^3}{|I_0|^3}
2^{-2m} \sum_T |I_T|.
\end{equation}
for all $I \subseteq I_0$.

Fix $I$.  The left-hand side of \eqref{quasi-outside} can be rewritten as
$$ \sum_\pv \sum_{\pv'} \| f_i \|_{P_i} \| f_i \|_{P'_i}
\int_{\R \backslash 2I_0} \Delta_{w_{P_i}}^* \phi_\pv(x)
\overline{\Delta_{w_{P'_i}}^* \phi_{\pv'}(x)} \tilde \chi_{I_0}^{-N}(x)\ dx,$$
where $\pv$, $\pv'$ are constrained by $I_\pv = I_{\pv'} = I$.

From the decay of $\phi_\pv$ and the kernel of $\Delta_{w_{P_i}}$, we may estimate the integral by $O( |I|^{N+1} |I_0|^{-N} )$.  By translating $w_{P_i}$ to be centered at the origin, and integrating by parts repeatedly, one can also obtain the bound of $|I| O(1 + |I| \dist(w_{P_i}, w_{P'_i}))^{-N}$.  Taking the geometric mean of these estimates, we can bound the left-hand side of \eqref{quasi-outside} by
$$ |I|^{N/2 + 1} |I_0|^{-N/2}
\sum_\pv \sum_{\pv'} \| f_i \|_{P_i} \| f_i \|_{P'_i} 
(1 + |I| \dist(w_{P_i}, w_{P'_i}))^{-N/2}.$$
By Schur's test (or Young's inequality), this is bounded by
$$ |I|^{N/2 + 1} |I_0|^{-N/2} \sum_\pv \| f_i \|_{P_i}^2.$$
Thus it is only left to show that
$$ \sum_\pv |I| \| f_i \|_{P_i}^2 \lesssim 
2^{-2m} \sum_T |I_T|.$$
But this follows from \eqref{sup-size} and the observation that each tree $T$ contributes at most $O(1)$ multi-tiles $\pv$ to the left-hand sum.

It thus remains to show that
\be{quasi-inside}
\| \sum_{\pv \in \bigcup_T T} \| f_i \|_{P_i} \Delta_{w_{P_i}}^* \phi_\pv
\|_2^2 \lesssim
2^{-2m} \sum_T |I_T|.
\end{equation}
We estimate the left-hand side of \eqref{quasi-inside} as
$$
\sum_{\pv, \pv' \in \bigcup_T T} \| f_i \|_{P_i} \|f_i \|_{P'_i}
|\langle \Delta_{w_{P_i}}^* \phi_\pv, \Delta_{w_{P'_i}}^* \phi_{\pv'}
\rangle|.$$
The inner product vanishes unless $w_{P_i}$ and $w_{P'_i}$ intersect;
by the nesting property of dyadic intervals this means that one of these intervals is a subset of the other.  By symmetry it suffices to consider the case $w_{P_i} \subseteq w_{P'_i}$.

One can easily verify that $\Delta_{w_{P_i}}^* \phi_\pv \lesssim \tilde \chi_{I_\pv}^{2N}$, and similarly with $\pv$ replaced by $\pv'$.  Thus we may estimate the inner product as
$$ |\langle \Delta_{w_{P_i}}^* \phi_\pv, \Delta_{w_{P'_i}}^* \phi_{\pv'}
\rangle| \lesssim |I_{\pv'}| (1 + \frac{\dist(I_{\pv'},I_\pv)}{|I_\pv|})^{-2N}.$$
To show \eqref{quasi-inside} it thus suffices to show that
\be{quasi-quasi}
\sum_{\pv, \pv' \in \bigcup_T T: w_{P_i} \subseteq w_{P'_i}} 
\| f_i \|_{P_i} \|f_i \|_{P'_i}
|I_{\pv'}| (1 + \frac{\dist(I_{\pv'},I_\pv)}{|I_\pv|})^{-2N}
\lesssim 2^{-2m} \sum_T |I_T|.
\end{equation}

Let us first deal with the portion of the sum where $|I_\pv| \sim |I_{\pv'}|$.
In this case we use the estimate 
$$\| f_i \|_{P_i} \| f_i \|_{P'_i} \lesssim \| f_i\|_{P_i}^2 + \|f_i \|_{P'_i}^2.$$
We treat the first term, as the second is similar.  For each $\pv$, the associated $\pv'$ have disjoint spatial intervals $I_{\pv'}$.  Thus one may compute the $\pv'$ summation, and estimate this contribution to \eqref{quasi-quasi} as
$$ \sum_{\pv \in \bigcup_T T} \| f_i \|_{P_i}^2 |I_\pv|.$$
But this is acceptable by \eqref{albedo}.

Now suppose $|I_\pv| \gg |I_{\pv'}|$.  By \eqref{sup-size} we may estimate the contribution to \eqref{quasi-quasi} by
$$ 2^{-2m}
\sum_T \sum_{\pv \in T} 
 \sum_{\pv' \in \bigcup_{T'} T': w_{P_i} \subseteq w_{P'_i}, |I_\pv| \gg |I_{\pv'}|} 
|I_{\pv'}| (1 + \frac{\dist(I_{\pv'},I_\pv)}{|I_\pv|})^{-2N}$$
From the assumptions on $\pv$ and $\pv'$ we see 
that $\pv'$ must belong to a tree other than $T$; since the trees are strongly $i$-disjoint we thus have $I_{\pv'} \cap I_T = \emptyset$, and that the $I_{\pv'}$ are disjoint.  We may thus estimate the contribution to \eqref{quasi-quasi} by
$$ 2^{-2m}
\sum_T \sum_{\pv \in T} 
\int_{\R \backslash I_T}
(1 + \frac{\dist(x,I_\pv)}{|I_\pv|})^{-2N}\ dx.$$
The integral bounded by
$$(1 + \frac{\dist(\R \backslash I_T,I_\pv)}{|I_\pv|})^{-3}.$$
Inserting this into the previous and computing the inner sum, we obtain \eqref{quasi-quasi} as desired.  This completes the proof of Lemma \ref{disjoint-lemma}.

\section{Proof of Lemma \ref{selection}}\label{select-sec}

Fix $i$, $m$.  The idea will be to remove trees $T$ from $\Pv$ one at a time until \eqref{remainder} is satisfied.

By refining the tree by a finite factor we may assume (using Lemma \ref{rank}) that for each dyadic interval $I$ there is at most one multi-tile $\pv \in T$ such that $I_\pv = I$.  We may assume that for any $\pv, \pv' \in \Pv$, $|I_\pv|/|I_\pv'|$ is an integer power of $2^{C_1}$, where $C_1$ is a large constant to be chosen shortly.  By Lemma \ref{rank} and a further refinement we can ensure that if $w_{P_i}$ is fixed, then
$w_{P_j}$ is also fixed for every $1 \leq j \leq n$.

Let $\Pv^*$ consist of those multi-tiles $\pv$ in $\Pv$ such that
$$ \| f_i \|_{P_i} \geq 2^{-m-2};$$
for these tiles we thus have
\be{heavy}
 \| f_i \|_{P_i} \sim 2^{-m}
\end{equation}
by \eqref{m-bound}.  We place a partial order $<$ on the multi-tiles in $\Pv^*$ by defining $\pv' < \pv$ if $P'_i < P_i$.  Let $\Pv^{**}$ be those tiles which are maximal with respect to this ordering.

By construction, the tiles $\{ P_i: \pv \in \Pv^{**}\}$ are disjoint.  From this, \eqref{heavy}, and Corollary \ref{lacey} we see that
\be{catch} \sum_{\pv \in \Pv^{**}} |I_\pv| \lesssim 2^{2m} \| f_i \tilde \chi_{I_0}^{N/2} \|_2^2 \lesssim 2^{2m} \min(1, \lambda_1 |I_0|).
\end{equation}
For each $\pv \in \Pv^{**}$ we associate the $i$-tree
$$ T = \{ \pv' \in \Pv^*: \pv' \leq \pv \}.$$
From \eqref{catch} we see that one can remove these trees $T$ from $\Pv$ and
place them into $\T$ while respecting \eqref{t-size}.  After removing these trees,
we have eliminated all elements of $\Pv^*$, so that we have
\be{sashay}
 \| f_i \|_{P_i} < 2^{-m-2}
\end{equation}
for all remaining multi-tiles $\pv$.

If $P$ is a tile, let $\xi_P$ denote the center of $w_P$.
If $P$ and $P'$ are tiles, we write $P' \lesssim^+ P$ if $P' \lesssim' P$ and
$\xi_{P'} > \xi_P$, and $P' \lesssim^-$ if $P' \lesssim' P$ and $\xi_{P'} < \xi_P$.  If $T$ is a tree, write $\xi_{T,i}$ for $\xi_{P_{T,i}}$.

We now perform the following algorithm.  We consider the set of all
trees $T$ in $\Pv$ such that 
\be{plus}
P_i \lesssim^+ P_{T,i} \hbox{ for all } \pv \in T
\end{equation}
and
\be{t-plus}
\sum_{\pv \in T} |I_\pv| \| f_j\|_{P_i}^2 
\geq 2^{-2m-5} |I_T|.
\end{equation}
If there are no trees obeying \eqref{plus} and \eqref{t-plus}, we terminate the algorithm.  Otherwise, we choose $T$ among all such trees so that $\xi_{T,i}$ is maximal, and that $T$ is maximal with respect to set inclusion.  Let $T'$ denote the $i$-tree
$$ T' = \{ \pv \in \Pv: P_i \leq P_{T,i} \}.$$
We remove both $T$ and $T'$ from $\Pv$, and add them to $\T$.  (These two trees are allowed to overlap).  Then one repeats the algorithm until we run out of trees obeying \eqref{plus} and \eqref{t-plus}.

Since $\Pv$ is finite, this algorithm terminates in a finite number of steps, producing trees $T_1, T'_1, T_2, T'_2, \ldots, T_M, T'_M$.  We claim that
the trees $T_1, \ldots, T_M$ produced in this manner are strongly disjoint.
It is clear from construction that $T_s \cap T_{s'} = \emptyset$ for all
$s \neq s'$; by our assumptions on the multi-tiles we thus see that
$P_i \neq P'_i$ for all $\pv \in T_s$, $\pv' \in T_{s'}$, $s \neq s'$.

Now suppose for contradiction that we had multi-tiles
$\pv \in T_s$, $\pv' \in T_{s'}$ such that $w_{P_i} \subsetneq w_{P'_i}$
and $I_{P'_i} \subseteq I_{T_s}$.  From our assumptions on the multi-tiles
we thus have $|w_{P'_i}| \geq 2^{C_1} |w_{P_i}|$.  Since $P_i \lesssim P_{T_s,i}$ and $P'_i \lesssim^+ P_{T_{s'},i}$, we thus see that
$\xi_{T_{s'},i} < \xi_{T_s,i}$ if $C_1$ is sufficiently large.  By our selection 
algorithm this implies that $s < s'$.

Also, since $|w_{P'_i}| \geq 2^{C_1} |w_{P_i}|$, $I_{P'_i} \subseteq I_{T_s}$, and $P_i \lesssim P_{T_s,i}$ we see that $P'_i \leq P_{T_s,i}$ if $C_1$ is
sufficiently large.  Since $s < s'$, this means that $\pv' \in T'_s$.  But
$T'_s$ and $T_{s'}$ are disjoint by construction, which is a contradiction.
Thus the trees $T_s$ are strongly disjoint.  From \eqref{m-bound} and
\eqref{t-plus} we see that these trees obey \eqref{size-lower}, and thus we have
$$ \sum_{s=1}^M |I_{T_s}| \lesssim 2^{2m} \min(1, \lambda_i |I_0|).$$
Since $T'_s$ has the same top as $T_s$, we may thus add all the $T_s$ and $T'_s$
to $\T$ while respecting \eqref{t-size}.

Now consider the set $\Pv$ of remaining multi-tiles.  We note that
\be{plus-left}
\sum_{\pv \in T: P_i \lesssim^+ P_{T,i}} |I_\pv| \| f_j\|_{P_i}^2 
< 2^{-2m-5} |I_T|
\end{equation}
for all trees $T$ in $\Pv$, since otherwise the portion of $T$ which obeyed \eqref{plus} would be eligible for selection by the above algorithm.

We now repeat the previous algorithm, but replace $\lesssim^+$ by $\lesssim^-$
and select the trees $T$ so that $\xi_{T,i}$ is \emph{minimized} rather
than maximized.  This yields a further collection of trees to add to $\T$ while still respecting \eqref{t-size}, and the remaining collection of tiles $\Pv$ has the property that
\be{plus-right}
\sum_{\pv \in T: P_i \lesssim^- P_{T,i}} |I_\pv| \| f_j\|_{P_i}^2 
< 2^{-2m-5} |I_T|
\end{equation}
for all trees $T$ in $\Pv$.  Combining \eqref{sashay}, \eqref{plus-left}, and \eqref{plus-right} we see that
$$ \size_i(T) \leq 2^{-m-1}$$
for all trees $T$ in $\Pv$, and we are done.

\section{Proof of Lemma \ref{lambda-lemma}}\label{lambda-sec}

Fix $1 \leq i \leq n$.  We may refine the collection $T$ of tiles as in the previous section.
 
Let $P$ be a tile.  Since the convolution kernel of $\Delta_{w_P}$ is rapidly decreasing for $|x| \gg |I_P|$, we see from the definition of $\|f\|_P$ that
$$ \| f\|_P \lesssim \frac{1}{|I_P|} \| f \tilde \chi_{I_P}^N \|_1.$$
From \eqref{F1} we thus have
$$ \| f_i \|_{P_i} \lesssim a_i \lambda_i $$
for all $\pv \in \Pv$.  In particular we have
$$ \| f_i \|_{T,i} \lesssim a_i \lambda_i.$$
for all trees $T$ in $\Pv$.

Let $B$ denote the best constant such that
\be{b-sir}
\sum_{\pv \in T: P_i \lesssim' P_{T,i}} |I_\pv| \| f_j\|_{P_i}^2
\leq B |I_T|
\end{equation}
for all trees $T$ in $\Pv$; to finish the proof of Lemma \ref{lambda-lemma} we must show that $B \lesssim a_i^2 \lambda_i^2$.

To achieve this we first need to prove an apparently weaker estimate.

\begin{lemma}\label{weak}  For any tree $T$ and function $f$, we have
$$
\| (\sum_{\pv \in T: P_i \lesssim' P_{T,i}} \| f\|_{P_i}^2 \chi_{I_\pv})^{1/2}
\|_{L^{1,\infty}} \lesssim \| f \tilde \chi_{I_T}^N \|_1.
$$
\end{lemma}

\begin{proof}  The expression in the norm is a variant of a Littlewood-Paley square function.  Thus, we shall use Calder\'on-Zygmund techniques to prove this estimate.

By frequency translation invariance we may assume that $w_{T,i}$ contains the origin.

Let us first assume that $f$ is supported outside of $2I_T$.  From Lemma \ref{triv} we have
\be{fp-decay}
\| f \|_{P_i} \lesssim \frac{|I_\pv|^{N-1}}{|I_T|^N} \| f \tilde \chi_{I_T}^N \|_1.
\end{equation}
Applying this estimate, we obtain
$$
\| (\sum_{\pv \in T: P_i \lesssim' P_{T,i}} \| f\|_{P_i}^2 \chi_{I_\pv})^{1/2} \|_2
\lesssim |I_T|^{-1/2} \| f \tilde \chi_{I_T}^N \|_1,$$
and the claim follows from H\"older.

It thus remains to show that
\be{wt2} \{ (\sum_{\pv \in T: P_i \lesssim' P_{T,i}} \| f\|_{P_i}^2 \chi_{I_\pv})^{1/2} 
\gtrsim \alpha \} \lesssim \alpha^{-1} \| f \|_1
\end{equation}
for all $\alpha > 0$.

Fix $\alpha$.  Perform a Calder\'on-Zygmund decomposition at level $\alpha$
$$ f = g + \sum_I b_I$$
where $\|g\|_2 \lesssim \alpha^{1/2} \|f\|_1^{1/2}$, the $I$ are intervals such that 
\be{i-sum}
\sum_I |I| \lesssim \alpha^{-1} \|f\|_1,
\end{equation} and the $b_I$ are supported on $I$ and satisfy $\int_I b_I \sim \alpha |I|$ and $\int b_I = 0$.

To control the contribution of $g$, it suffices from Chebyshev to verify the $L^2$ bound
$$ \| (\sum_{\pv \in T: P_i \lesssim' P_{T,i}} \| g\|_{P_i}^2 \chi_{I_\pv})^{1/2} \|_2
\lesssim \|g\|_2.$$
The left-hand side of this is
\be{lacun} (\sum_{\pv \in T: P_i \lesssim' P_{T,i}} |I_\pv| \| g\|_{P_i}^2)^{1/2}.
\end{equation}
However, from H\"older and the definition of $\|g\|_{P_i}$ we have
$$
|I_\pv| \| g\|_{P_i}^2 \lesssim \| \Delta_{w_{P_i}}(g) \tilde \chi_{I_\pv}^{-1} \|_2^2.$$
Thus we may bound \eqref{lacun} by
$$ (\sum_w \| \Delta_w(g) \|_2^2)^{1/2},$$
where $w$ ranges over the set $\{ w_{P_i}: \pv \in T \}$.  But the desired bound of $\|g\|_2$ then follows from Plancherel and the lacunary nature of the $w$.

To deal with the $b_I$, it suffices from the triangle inequality, Chebyshev, and
to show that
$$ \| (\sum_{\pv \in T: P_i \lesssim' P_{T,i}} \| b_I\|_{P_i}^2 \chi_{I_\pv})^{1/2} \|_{L^1(\R \backslash 2I)} \lesssim \alpha |I|$$
for all $I$.  In fact we prove the stronger
\be{ok} \| \sum_{\pv \in T: P_i \lesssim' P_{T,i}} \| b_I\|_{P_i} \chi_{I_\pv} \|_{L^1(\R \backslash 2I)} \lesssim \alpha |I|.
\end{equation}

Fix $I$.  We may restrict the summation to those $\pv$ such that
$I_\pv \not \subseteq 2I$.  

From Lemma \ref{triv} we have
$$ \| b_I \|_{P_i} \lesssim \alpha \frac{|I|}{|I_\pv|}
(1 + \frac{\dist(I_\pv,I)}{|I_\pv|})^{-N};$$
in particular, from the hypothesis $I_\pv \not \subseteq 2I$ we have
$$ \| b_I \|_{P_i} \lesssim \alpha \frac{|I_\pv|^{N-1}}{|I|^{N-1}}.$$
Also, by playing off the moment condition on $b_I$ against the smoothness of
$\Delta_{w_{P_i}}$, we have
$$ \| b_I \|_{P_i} \lesssim \alpha \frac{|I|^2}{|I_\pv|^2}.$$
Combining all these estimates, we obtain
$$ \| b_I \|_{P_i} \lesssim \alpha \frac{|I|}{|I_\pv|}
(1 + \frac{\dist(I_\pv,I)}{|I_\pv|})^{-N/2}
\min( \frac{|I_\pv|}{|I|}, \frac{|I|}{|I_\pv|})^{1/2}.$$
Inserting this into \eqref{ok} we obtain the result.
\end{proof}

To bootstrap Lemma \ref{weak} to Lemma \ref{lambda-lemma} we shall employ a variant of arguments used to prove the John-Nirenberg inequality.

By construction of $B$, there exists a tree $T$ such that
\be{max}
\sum_{\pv \in T: P_i \lesssim' P_{T,i}} |I_\pv| \| f_j\|_{P_i}^2
= B |I_T| 
\end{equation}

Fix this tree.  From Lemma \ref{weak} and \eqref{F1} we have
$$ \| (\sum_{\pv \in T: P_i \lesssim' P_{T,i}} \| f\|_{P_i}^2 \chi_{I_\pv})^{1/2}
\|_{L^{1,\infty}} \lesssim |I_T| a_i \lambda_i.$$
We thus have $|E| \leq \frac{1}{2} |I_T|$,
where
$$ E = \{ \sum_{\pv \in T: P_i \lesssim' P_{T,i}} \| f\|_{P_i}^2 \chi_{I_\pv}
\geq C a_i^2 \lambda_i^2 \}$$
and $C$ is a sufficiently large constant.

From the nesting properties of dyadic intervals we see that there must exist  a subset $T^*$ of $T$ such that the intervals $\{ I_\pv: \pv \in T^* \}$ form
a partition of $E$.  In particular we have
\be{ch}
\sum_{\pv \in T^*} |I_\pv| \leq \frac{1}{2} |I_T|.
\end{equation}

If $\pv \in \Pv$ is such that $I_\pv \not \subseteq E$, then we must have $P_i \lesssim' P'_i$ for some $\pv' \in T^*$, if $C_1$ is chosen sufficiently large.
We can thus decompose the left-hand side of \eqref{max} as
$$
\| \sum_{\pv \in T: P_i \lesssim' P_{T,i}, I_\pv \not \subseteq E} \| f_j\|_{P_i}^2 \chi_{I_\pv} \|_1
+ \sum_{\pv' \in T^*} \sum_{\pv \in T: P_i \lesssim' P_{T,i}, I_\pv \subseteq E} |I_\pv| \| f_j\|_{P_i}^2.$$
Consider the former term.  From the definition of $E$ and the nesting properties of dyadic intervals we see that the expression in the norm is $O(a_i^2 \lambda_i^2)$.  Thus the former term is $O(|I_T| a_i^2 \lambda_i^2)$.

Now consider the latter summation.  For each $\pv' \in T^*$ the inner sum is
$O(B |I_{\pv'}|)$ from \eqref{b-sir}.  Inserting these estimates back into \eqref{max} and using \eqref{ch} we obtain
$$ B |I_T| \lesssim |I_T| a_i^2 \lambda_i^2
+ \frac{1}{2} B |I_T|,$$
and the claim follows.
This concludes the proof of Lemma \ref{lambda-lemma}.

\section{Conclusion of the $k=1$ case}\label{conclusion}

We now prove \eqref{disc}.  We first observe from iterating Lemma \ref{selection} and using Lemma \ref{lambda-lemma} that

\begin{corollary}  Let $1 \leq i \leq n$.  Then there exists a partition
$$ \Pv = \bigcup_{m: 2^{-m} \lesssim \lambda_i a_i} \Pv^{m,i}$$
where one has \eqref{size-lower} for all trees $T$ in $\Pv^{m,i}$, and
such that $\Pv^{m,i}$ can be covered as
\be{cover}
\Pv^{m,i} = \bigcup_{T \in \T^{m,i}} T
\end{equation}
where $\T^{m,i}$ is a collection of trees such that
\be{tm-size}
\sum_{T \in \T^{m,i}} |I_T| \lesssim 2^{2m} \min(1, \lambda_i |I_0|).
\end{equation}
\end{corollary}

Write the left-hand side of \eqref{disc} as
$$
\sum_{m_1, \ldots, m_n} \sum_{\pv \in \Pv^{m_1,1} \cap \ldots \cap
\Pv^{m_n,n}} |I_\pv| \prod_{i=1}^n \|f_i\|_{P_i}
$$
where we implicitly assume
\be{moo}
2^{-m_i} \leq \lambda_i a_i.
\end{equation}
By symmetry we may restrict the summation to the case
\be{monotone}
m_1 \leq m_2 \leq \ldots \leq m_n.
\end{equation}
We then estimate the sum by
\be{cor-split}
\sum_{m_1 \leq \ldots \leq m_n} \sum_{T \in \T^{m_1,1}} 
\sum_{\pv \in T'} |I_\pv| \prod_{i=1}^n \|f_i\|_{P_i}
\end{equation}
where $T' = T'(T, m_2, \ldots, m_n)$ denotes the tree
$$ T' = T \cap \Pv^{m_2,2} \cap \ldots \cap \Pv^{m_n,n}.$$

By Lemma \ref{Cauchy-Schwarz} we may estimate \eqref{cor-split} by
\be{cor-bound}
\sum_{m_1 \leq \ldots \leq m_n} \sum_{T \in \T^{m_1,1}} 
|I_T| \sup_{1 \leq i_1 < i_2 \leq n}
\size_{i_1}(T') \size_{i_2}(T') \prod_{i \neq i_1,i_2} \| f_i \|_{T',i}.
\end{equation}

From \eqref{size-lower} we have
$$ \size_i(T') \lesssim 2^{-m_i},$$
which implies with \eqref{a-def} that
$$ \| f_i \|_{T',i} \lesssim \min(2^{-m_i},A).$$
Thus we may estimate \eqref{cor-bound} by
$$
\sum_{m_1 \leq \ldots \leq m_n} \sum_{T \in \T^{m_1,1}} 
|I_T| \sup_{1 \leq i_1 < i_2 \leq n}
2^{-m_{i_1}} 2^{-m_{i_2}} \prod_{i \neq i_1,i_2} \min(2^{-m_i},A).
$$
It is clear that the supremum is attained when $i_1 = 1$, $i_2 = 2$.
Applying \eqref{tm-size} we can thus estimate the previous by
$$
\min(1,\lambda_1 |I_0|) \sum_{m_1 \leq \ldots \leq m_n} 2^{2m_1} 2^{-m_1} 2^{-m_2} \prod_{2 < i \leq n} \min(2^{-m_i},A).$$
Clearly we have the estimate
$$ \min(1,\lambda_1 |I_0|) \leq \min(1,|I_0|) \prod_{i=1}^n (1 + \lambda_i).$$
To show \eqref{disc}, it thus suffices to show
\be{disc-a}
\sum_{m_1 \leq \ldots \leq m_n} 2^{2m_1} 2^{-m_1} 2^{-m_2} \prod_{2 < i \leq n} \min(2^{-m_i},A)
\lesssim
A^{n-2-\theta_1 - \ldots - \theta_n} 
\prod_{i=1}^n (a_i \lambda_i)^{\theta_i}.
\end{equation}

We first consider the case when \eqref{theta-ineq} holds with equality (i.e.
\eqref{theta-bound} holds).  In this case we need only show that
\be{soosh}
\sum_{m_1 \leq \ldots \leq m_n} 2^{2m_1} \prod_{i=1}^n 2^{-m_i}
\lesssim
\prod_{i=1}^n (a_i \lambda_i)^{\theta_i}.
\end{equation}
From \eqref{theta-ineq} we may write
\be{pest} 2^{2m_1} = \prod_{i=1}^n 2^{(1-\theta_i) m_1} \leq
\prod_{i=1}^n 2^{(1-\theta_i) m_i}
\end{equation}
by \eqref{theta} and \eqref{monotone}.  Thus \eqref{soosh} reduces to
\be{cocky}
\sum_{m_1, \ldots, m_n} \prod_{i=1}^n 2^{-\theta_i m_i}
\lesssim
\prod_{i=1}^n (a_i \lambda_i)^{\theta_i}.
\end{equation}
But this follows from \eqref{moo} and \eqref{theta}.

Now suppose that \eqref{theta-ineq} holds with strict inequality.  We may then find $\theta'_i$ satisfying \eqref{theta} and \eqref{theta-bound} such that $\theta_i = \theta'_i$ for $i=1,2$ and $\theta'_i > \theta_i$ for $i > 2$;
note how one needs \eqref{theta} and \eqref{theta-ineq} for $k=1$ to ensure that $\theta'_i$ exists. 
Using the estimate
$$ \min(2^{-m_i},A) \leq A^{\theta'_i - \theta_i} 2^{-m_i(1 - \theta'_i + \theta_i)}$$
and canceling the $A$ factors, we reduce to
$$
\sum_{m_1 \leq \ldots \leq m_n} 2^{2m_1} \prod_{i=1}^n 2^{-m_i(1 - \theta'_i + \theta_i)}
\lesssim
\prod_{i=1}^n (a_i \lambda_i)^{\theta_i}.
$$
Applying \eqref{pest} with the $\theta_i$ replaced by $\theta'_i$, we reduce to
\eqref{cocky} as before.  Thus in either case \eqref{disc} is proven.

\section{The induction on $k$}\label{induct-sec}

We have just proven Theorem \ref{discretized} when $k=1$.  Now suppose inductively that $k > 1$, and the claim has already been proven for $k-1$.

We need to show \eqref{disc}.  By symmetry it suffices to consider those tiles $\pv$ for which
\be{mono-tile}
\| f_1 \|_{P_1} \geq \| f_2 \|_{P_2} \geq \ldots \geq \| f_n \|_{P_n};
\end{equation}
we shall implicitly assume this in the sequel.

From Lemma \ref{triv}, \eqref{F1}, and \eqref{a-def} we have
$$ \|f_1\|_{P_1} \lesssim \min(a_1 \lambda_1,A).$$
Thus \eqref{disc} reduces to showing that
\be{disc-m} \sum_{m: 2^{-m} \lesssim \min(a_1 \lambda_1,A)}
\sum_{\pv \in \Pv^m}
|I_\pv| \prod_{i=1}^n \|f_i\|_{P_i}
\lesssim A^{n-2k-\theta_1 - \ldots - \theta_n} \min(1, |I_0|)
\prod_{i=1}^n (\lambda_i a_i)^{\theta_i} (1 + \lambda_i)
\end{equation}
where
$$ \Pv^m = \{ \pv \in \Pv: \| f_1 \|_{P_1} \sim 2^{-m} \}.$$

Fix $m$.  We order the multi-tiles in $\Pv^m$ by setting $\pv' < \pv$
if $P'_1 < P_1$.  Let $\Pv^{m,*}$ be the tiles in $\Pv^m$ which are maximal with respect to this ordering.  By applying Corollary \ref{lacey} as in the proof of
Lemma \ref{selection}, we see that
\be{asok}
 \sum_{\pv \in \Pv^{m,*}} |I_\pv| \lesssim 2^{2m} \min(1, \lambda_1 |I_0|).
\end{equation}

We may estimate the left-hand side of \eqref{disc-m} as
\be{step}
\sum_{m: 2^{-m} \lesssim \min(a_1 \lambda_1,A)}
\sum_{\pv' \in  \Pv^{m,*}} \sum_{\pv \in \Pv^m: P_1 \leq P'_1}
|I_\pv| 2^{-m} \prod_{i=2}^n \|f_i\|_{P_i}.
\end{equation}
For fixed $\pv'$, the collection of multi-tiles $\{ \pv \in \Pv^m: P_1 \leq P'_1 \}$ satisfies the conditions of Lemma \ref{rank} with $(n,k)$ replaced by $(n-1,k-1)$, if we forget the first tile $P_1$ from each multi-tile $\pv$.  Thus we may apply the induction hypothesis, with $I_0$ replaced by $I_{\pv'}$ and $A$
estimated by $2^{-m}$ (thanks to \eqref{mono-tile}), and estimate
\eqref{step} by
$$
\sum_{m: 2^{-m} \lesssim \min(a_1 \lambda_1,A)} 
\sum_{\pv' \in  \Pv^{m,*}} 2^{-m} \min(1,|I_{\pv'}|)
2^{-m((n-1)-2(k-1)-\theta_2 - \ldots - \theta_n)}
\prod_{i=2}^n (\lambda_i a_i)^{\theta_i} (1+\lambda_i).$$
Estimating $\min(1,|I_{\pv'}|)$ by $|I_{\pv'}|$ and applying \eqref{asok}, and
then gathering the powers of $2^m$, this can be estimated by
$$
\sum_{m: 2^{-m} \lesssim \min(a_1 \lambda_1,A)} 
2^{-m(n-2k- \theta_2 - \ldots - \theta_n)}
\min(1, \lambda_1|I_0|)
\prod_{i=2}^n (\lambda_i a_i)^{\theta_i} (1 + \lambda_i).$$
Evaluating the $m$ summation and applying the elementary inequalities
$$ 
\min(a_1 \lambda_1,A)^{n-2k- \theta_2 - \ldots - \theta_n}
\lesssim A^{n-2k-\theta_2 - \ldots - \theta_n} (a_1 \lambda_1)^{\theta_1}
$$
(which follows from \eqref{theta} and \eqref{theta-ineq}) and
$$ \min(1, \lambda_1 |I_0|) \leq \min(1, |I_0|) (1 + \lambda_1)$$
we see that \eqref{disc-m} follows.  This concludes the proof of Theorem \ref{discretized}, and thus Theorem \ref{main}, for general $k$.

\section{Remarks}\label{remarks}

Let $A$ denote the Wiener algebra, that is the space of functions whose Fourier transform is in $L^1$.  The purpose of this section is to extend Theorem \ref{main} slightly to

\begin{theorem} Let $0 \leq s < n-1$, and let $\subspace$ be a subspace of $\hyperplane$ of dimension $k$ where
$$
0 \leq k-s < (n-s)/2. 
$$
Assume that $\subspace$ is non-degenerate in the sense of Theorem \ref{main}, and that $m$ satisfies \eqref{symb-gamma}.  Then one has
\be{holder-wiener}
T: L^{p_1} \times \ldots \times L^{p_{n-s-1}} \times A \times \ldots \times A
\to L^{p'_{n-s}}
\end{equation}
$1 < p_i \le \infty$ for $i=1, \ldots, n-s-1$, 
$$ \frac{1}{p_1} + \ldots + \frac{1}{p_{n-s}} = 1,$$
and
$$
\frac{1}{p_{i_1}} + \ldots + \frac{1}{p_{i_r}} < \frac{(n-s)-2(k-s)+r}{2}
$$
for all $1 \leq i_1 < \ldots < i_r \leq n-s$ and $1 \leq r \leq n-s$.
\end{theorem}

Thus, for instance, the trilinear Hilbert transform maps $L^p \times L^q \times A$ to $L^r$ whenever $1 < p,q \le \infty$, $1/p + 1/q = 1/r$, and 
$2/3<r <\infty$.

\begin{proof}  Let $g_1, \ldots, g_s$ be elements in the unit ball of $A$, and let $T_g$ denote the $n-s-1$-linear operator
$$ T_g(f_1, \ldots, f_{n-s-1}) = T(f_1, \ldots, f_{n-s-1},g_1, \ldots, g_s),$$
and let $\Lambda_g$ be the associated $n-s$-form as in \eqref{lambda-t}.
We need to show that
$$ T_g: L^{p_1} \times \ldots \times L^{p_{n-s-1}} \to L^{p'_{n-s}}.$$
By the reductions in Section \ref{interp-sec} it suffices to show that $\Lambda_g$ is of restricted type $p$ for all exponent $n-s$-tuples $p$
such that
$1 < p_i < 2$
for all indices $i$ which are not equal to the bad index $j$ of $p$, and 
$$
(k-s) - (n-s) + 2 > \frac{1}{p_j} > (k-s) - (n-s) + \frac{3}{2}.
$$

Fix $p$; by symmetry we may assume that $p$ has bad index $n-s$.  Let $E_1, \ldots, E_{n-s}$ be sets of finite measure.  We have to find a major subset $E'_{n-s}$ of $E_{n-s}$ such that
\be{lag} |\Lambda_g(F_1, \ldots, F_{n-s})| \lesssim |E|^{1/p}
\end{equation}
for all $F_i \in X(E'_i)$.  By scaling we may take $|E_{n-s}| = 1$.

We choose $E'_{n-s}$ to be the set defined by \eqref{en-def}, with $n$ replaced by $n-s$ throughout.   Since \eqref{lag} is sub-additive in $g$, and the unit ball of $A$ is the convex hull of the plane waves, we may assume that each $g_j$ is a plane wave $g_j(x) = e^{2\pi i x \xi_j}$ for some constants $\xi_j$.  By modulating the $F_i$ suitably, and translating the symbol $m$ by a direction in $\subspace$, one may set $\xi_j = 0$.  The functions $g$ are now completely harmless, and the claim follows from Theorem \ref{rwt-bound} with $n$ replaced by $n-s$.
\end{proof}

\end{document}